\documentclass[10pt,onecolumn]{article}
\usepackage{multicol}
\usepackage[greek,english]{babel}
\usepackage[vcentering,dvips]{geometry}
\usepackage{amsfonts}
\usepackage{amssymb}
\usepackage{eucal}
\usepackage{amsxtra}
\usepackage{color}
\usepackage{setspace}
\usepackage{times}
\usepackage{titlesec}
\usepackage{sectsty}
\allsectionsfont{\large}
\usepackage{url}
\usepackage{doi}
\usepackage{hyperref}
\hypersetup{
    colorlinks,%
    citecolor=black,%
    filecolor=black,%
    linkcolor=black,%
    urlcolor=black
}

\addtocounter{section}{0} \numberwithin{equation}{section}
\newtheorem{theorem}{Theorem}[section]
\newtheorem{proposition}[theorem]{Proposition}
\newtheorem{lemma}[theorem]{Lemma}

\newcommand{\qed}{{$\hfill \Box$}}

\onehalfspace

\hypersetup{urlcolor=blue }

\begin{document}
\pagestyle{plain}

\title{\textbf{Real hypersurfaces in $\mathbb{C}P^{2}$ and
$\mathbb{C}H^{2}$ whose structure Jacobi operator is Lie $\mathbb{D}$-parallel}}

\author{ \textbf{\normalsize{Konstantina Panagiotidou and Philippos J. Xenos}}\\
\small \emph{Mathematics Division-School of Technology, Aristotle University of Thessaloniki, Thessaloniki 54124, Greece}\\
\small \emph{E-mail: kapanagi@gen.auth.gr, fxenos@gen.auth.gr}}

\date{}

\maketitle

\begin{flushleft}
\small {\textsc{Abstract}.
In [3],[7] and [8] results concerning the parallelness of the Lie derivative of the structure Jacobi operator of a real hypersurface with respect to $\xi$ and to any vector field  $X$ were obtained in both complex projective space and complex hyperbolic space. In the present paper, we study the parallelness of the Lie derivative of the structure Jacobi operator of a real hypersurface with respect to vector field $X$ $\epsilon$ $\mathbb{D}$ in $\mathbb{C}P^{2}$ and $\mathbb{C}H^{2}$. More precisely, we prove that such real hypersurfaces do not exist.}
\end{flushleft}
\begin{flushleft}
\small{\emph{Keywords}: Real hypersurface, Lie $\mathbb{D}$-parallelness,
Structure Jacobi operator, Complex projective space, Complex hyperbolic space.\\}
\end{flushleft}
\begin{flushleft}
\end{flushleft}
\footnote{\thanks{The first author is granted by the Foundation Alexandros S. Onasis. Grant Nr: G ZF 044/2009-2010.}}

\section{Introduction}
\hspace{15pt}A complex n-dimensional Kaehler manifold of constant holomorphic
 sectional curvature c is called a complex space form, which is
 denoted by $M_{n}(c)$. A complete and simply connected
 complex space form is complex analytically isometric to a complex
 projective space $\mathbb{C}P^{n}$, a complex Euclidean space
 $\mathbb{C}^{n}$ or a complex hyperbolic space $\mathbb{C}H^{n}$
 if $c>0, c=0$ or $c<0$ respectively.

Let M be a real
 hypersurface in  a complex space form $M_{n}(c)$, $c\neq0$. Then an almost contact metric
 structure $(\varphi,\xi,\eta,g)$ can be defined on M induced from the Kaehler metric and complex structure J on $M_{n}(c)$.
  The structure vector field $\xi$ is called principal if $A\xi=\alpha\xi$, where A is the
 shape operator of M and $\alpha=\eta(A\xi)$ is a smooth function. A real hypersurface is said to
 be a \textit{Hopf hypersurface} if $\xi$ is principal.

The classification problem of real hypersurfaces in complex space forms is of great importance in Differential Geometry. The study of this was initiated by Takagi (see \cite{T}), who
classified homogeneous real hypersurfaces in $\mathbb{C}P^{n}$ and
showed that they could be divided into six types, which are said
to be of type $A_{1}$, $A_{2}$, $B$, $C$, $D$ and $E$. Berndt (see \cite{Ber})
classified homogeneous real hypersurfaces in $\mathbb{C}H^{n}$
with constant principal curvatures.

The Jacobi operator with respect to \emph{X} on M is defined by
 $R(\cdot, X)X$, where R is the Riemmanian curvature of M. For
 $X=\xi$ the Jacobi operator is called structure Jacobi operator
 and is denoted by $l=R(\cdot, \xi)\xi$. It has a
 fundamental role in almost contact manifolds.
 Many differential geometers have studied real hypersurfaces in terms of the structure Jacobi operator.

 The study of real hypersurfaces whose structure Jacobi operator satisfies conditions concerned to the parallelness of it is a problem of great importance. In \cite{OPS} the nonexistence
of real hypersurfaces in nonflat complex space form with parallel
structure Jacobi operator ($\nabla l=0$) was proved. In \cite{PSaSuh} a
weaker condition ($\mathbb{D}$-parallelness, where $\mathbb{D}=ker(\eta)$), that is
$\nabla_{X}l=0$ for any vector field $X$ orthogonal to $\xi$, was
studied and it was proved the nonexistence of such hypersurfaces
in case of $\mathbb{C}P^{n}$ ($n\geq3$). The $\xi$-parallelness of
structure Jacobi operator in combination with other conditions was
another problem that was studied by many authors such as Ki, Perez, Santos, Suh (\cite{KPSaSuh}).

The Lie derivative of the structure Jacobi operator is another condition that has been studied extensively. More precisely, in \cite{PSa} proved the non-existence of real hypersurfaces in $\mathbb{C}P^{n}$, ($n\geq3$), whose Lie derivative of the structure Jacobi operator with respect to any vector field X vanishes (i.e. $L_{X}l=0$). On the other hand, real hypersurfaces in $\mathbb{C}P^{n}$, ($n\geq3$), whose Lie derivative of the structure Jacobi operator with respect to $\xi$ vanishes (i.e. $L_{\xi}l=0$, Lie $\xi$-parallel) are classified (see \cite{PSaSuh1}). Ivey and Ryan in \cite{IR} extend some of the above results in $\mathbb{C}P^{2}$ and $\mathbb{C}H^{2}$. More precisely, they proved that in $\mathbb{C}P^{2}$ and $\mathbb{C}H^{2}$ there exist no real hypersurfaces satisfying condition $L_{X}l=0$, for any vector field $X$, but real hypersurfaces satisfying condition $L_{\xi}l=0$ exist and they classified them. Additional, they proved that there exist no real hypersurfaces in $\mathbb{C}P^{n}$ or $\mathbb{C}H^{n}$, ($n\geq3$), satisfying condition $L_{X}l=0$, for any vector field $X$.

Following the notion of \cite{PSaSuh1}, the structure Jacobi operator is said to be \emph{Lie $\mathbb{D}$-parallel}, when the Lie derivative of it with respect to any vector field X $\epsilon$ $\mathbb{D}$ vanishes. So the following question raises naturally:

"\emph{Do there exist real hypersurfaces in non-flat $M_{2}(c)$ with Lie $\mathbb{D}$-parallel structure Jacobi operator?}"

In this paper, we study the above question in $\mathbb{C}P^{2}$ and
$\mathbb{C}H^{2}$.  The condition of \emph{Lie $\mathbb{D}$-parallel structure Jacobi operator}, i.e. $L_{X}l=0$ with X $\epsilon$ $\mathbb{D}$, implies:
\begin{eqnarray}
\nabla_{X}(lY)+l\nabla_{Y}X=\nabla_{lY}X+l\nabla_{X}Y,
\end{eqnarray}
where Y $\epsilon$ $TM$.\\

We prove the following:\\

\begin{pro}
There exist no real hypersurfaces in $\mathbb{C}P^{2}$ and
$\mathbb{C}H^{2}$ equipped with Lie $\mathbb{D}$-parallel structure Jacobi operator.
\end{pro}
\section{Preliminaries}

\hspace{15pt}Throughout this paper all manifolds, vector fields e.t.c are assumed to be of class $C^{\infty}$ and all manifolds are assumed to be connected. Let M be a connected real hypersurface immersed in a nonflat
complex space form $(M_{n}(c),G)$ with almost complex structure J
of constant holomorphic sectional curvature c. Let N be a locally defined unit
normal vector field on M and $\xi=-JN$. For a vector field X
tangent to M we can write $JX=\varphi (X)+\eta(X)N$, where
$\varphi X$ and $\eta(X)N$ are the tangential and the normal
component of $JX$ respectively. The Riemannian connection
$\overline{\nabla}$ in $M_{n}(c)$ and $\nabla$ in M are related
for any vector fields X,Y on M:
$$\overline{\nabla}_{Y}X=\nabla_{Y}X+g(AY,X)N$$
$$\overline{\nabla}_{X}N=-AX$$
where g is the Riemannian metric on M induced from G of $M_{n}(c)$
and A is the shape operator of M in $M_{n}(c)$. M has an almost
contact metric structure $(\varphi,\xi,\eta)$ induced from J on
$M_{n}(c)$ where $\varphi$ is a (1,1) tensor field and $\eta$ a
1-form on M such that ([2])
$$g(\varphi X,Y)=G(JX,Y),\hspace{20pt}\eta(X)=g(X,\xi)=G(JX,N).$$
Then we have
\begin{eqnarray}
\varphi^{2}X=-X+\eta(X)\xi,\hspace{20pt}
\eta\circ\varphi=0,\hspace{20pt} \varphi\xi=0,\hspace{20pt}
\eta(\xi)=1
\end{eqnarray}
\begin{eqnarray}\hspace{20pt}
g(\varphi X,\varphi
Y)=g(X,Y)-\eta(X)\eta(Y),\hspace{10pt}g(X,\varphi Y)=-g(\varphi
X,Y)
\end{eqnarray}
\begin{eqnarray}
\nabla_{X}\xi=\varphi
AX,\hspace{20pt}(\nabla_{X}\varphi)Y=\eta(Y)AX-g(AX,Y)\xi
\end{eqnarray}
    Since the ambient space is of constant holomorphic sectional
curvature c, the equations of Gauss and Codazzi for any vector
fields X,Y,Z on M are respectively given by
\begin{eqnarray}
R(X,Y)Z=\frac{c}{4}[g(Y,Z)X-g(X,Z)Y+g(\varphi Y ,Z)\varphi
X\end{eqnarray} $$-g(\varphi X,Z)\varphi Y-2g(\varphi X,Y)\varphi
Z]+g(AY,Z)AX-g(AX,Z)AY$$
\begin{eqnarray}
\hspace{10pt}
(\nabla_{X}A)Y-(\nabla_{Y}A)X=\frac{c}{4}[\eta(X)\varphi
Y-\eta(Y)\varphi X-2g(\varphi X,Y)\xi]
\end{eqnarray}
where R denotes the Riemannian curvature tensor on M.\\
    Relation (2.4) implies that the structure Jacobi operator $l$ is given by:
\begin{eqnarray}
lX=\frac{c}{4}[X-\eta(X)\xi]+\alpha AX-\eta(AX)A\xi
\end{eqnarray}

    For every point $P\;\;\epsilon\;\; M$, the tangent space
$T_{P}M$ can be decomposed as following:
$$T_{P}M=span\{\xi\}\oplus ker(\eta),$$
where $ker(\eta)=\{X\;\;\epsilon\;\; T_{P}M:\eta(X)=0\}$.
 Due to the above decomposition,the vector field $A\xi$ can be written:
 \begin{eqnarray}
 A\xi=\alpha\xi+\beta U,
 \end{eqnarray}
 where $\beta=|\varphi\nabla_{\xi}\xi|$ and
 $U=-\frac{1}{\beta}\varphi\nabla_{\xi}\xi\;\epsilon\;ker(\eta)$, provided
 that $\beta\neq0$.

\section{Some Previous Results}
Let M be a non-Hopf hypersurface in $\mathbb{C}P^{2}$ or $\mathbb{C}H^{2}$ (i.e. $M_{2}(c)$, $c\neq0$). Then the following relations holds on every three-dimensional real hypersurface in $M_{2}(c)$.
\begin{lemma}
Let M be a real hypersurface in $M_{2}(c)$. Then the following
relations hold on M:
\begin{eqnarray}
\hspace{-70pt}AU=\gamma U+\delta\varphi U+\beta\xi,\hspace{20pt}
A\varphi U=\delta U+\mu\varphi U,
\end{eqnarray}
\begin{eqnarray}
\nabla_{U}\xi=-\delta U+\gamma\varphi U,\hspace{20pt}
\nabla_{\varphi U}\xi=-\mu U+\delta\varphi U,\hspace{20pt}
\nabla_{\xi}\xi=\beta\varphi U,
\end{eqnarray}
\begin{eqnarray}
\nabla_{U}U=\kappa_{1}\varphi U+\delta\xi,\hspace{20pt}
\nabla_{\varphi U}U=\kappa_{2}\varphi U+\mu\xi,\hspace{20pt}
\nabla_{\xi}U=\kappa_{3}\varphi U,
\end{eqnarray}
\begin{equation}
\nabla_{U}\varphi U=-\kappa_{1}U-\gamma\xi,\hspace{5pt}
\nabla_{\varphi U}\varphi U=-\kappa_{2}U-\delta\xi,\hspace{5pt}
\nabla_{\xi}\varphi U=-\kappa_{3}U-\beta\xi,
\end{equation}
where $\gamma,\delta,\mu,\kappa_{1},\kappa_{2},\kappa_{3}$ are
smooth functions on M.
\end{lemma}
\textbf{Proof:} Let $\{U,\varphi U,\xi\}$ be an orthonormal basis of M. Then we
have:
$$\hspace{10pt}AU=\gamma U+\delta\varphi U+\beta\xi\hspace{30pt}A\varphi U=\delta U+\mu\varphi U,$$ where $\gamma,\delta,\mu$ are smooth functions, since $g(AU,\xi)=g(U,A\xi)=\beta$ and
$g(A\varphi U,\xi)=g(\varphi U,A\xi)=0$.

The first relation of (2.3), because of (2.6) and (3.1), for $X=U$,
$X=\varphi U$ and $X=\xi$ implies (3.2), owing to (2.7).

 From the well known
relation: $Xg(Y,Z)=g(\nabla_{X}Y,Z)+g(Y,\nabla_{X}Z)$ for
$X,Y,Z$ $\epsilon$ $\{\xi,U,\varphi U\}$ we obtain (3.3) and
(3.4), where $\kappa_{1}, \kappa_{2}$ and $\kappa_{3}$ are smooth
functions.\qed
\\

Because of Lemma 3.1 the Codazzi equation implies:
\begin{eqnarray}
U\beta-\xi\gamma&=&\alpha\delta-2\delta\kappa_{3}\\
\xi\delta&=&\alpha\gamma+\beta\kappa_{1}+\delta^{2}+\mu\kappa_{3}+\frac{c}{4}-\gamma\mu-\gamma\kappa_{3}-\beta^{2}\\
U\alpha-\xi\beta&=&-3\beta\delta\\
\xi\mu&=&\alpha\delta+\beta\kappa_{2}-2\delta\kappa_{3}\\
(\varphi U)\alpha&=&\alpha\beta+\beta\kappa_{3}-3\beta\mu\\
(\varphi U)\beta&=&\alpha\gamma+\beta\kappa_{1}+2\delta^{2}+\frac{c}{2}-2\gamma\mu+\alpha\mu\\
U\delta-(\varphi U)\gamma&=&\mu\kappa_{1}-\kappa_{1}\gamma-\beta\gamma-2\delta\kappa_{2}-2\beta\mu\\
U\mu-(\varphi U)\delta&=&\gamma\kappa_{2}+\beta\delta-\kappa_{2}\mu-2\delta\kappa_{1}
\end{eqnarray}

We recall the following Proposition (\cite{IR}):
\begin{proposition}
There does not exist real non-flat hypersurface in $M_{2}(c)$, whose structure Jacobi operator vanishes.
\end{proposition}

\section{Auxiliary Relations}
If M is a real non-flat hypersurface in $\mathbb{C}P^{2}$ or
$\mathbb{C}H^{2}$ (i.e. $M_{2}(c)$, $c\neq0$), we  consider the open subset $\mathcal{W}$ of points
$P$ $\epsilon$ $M$, such that there exists a neighborhood of every $P$,
where $\beta=0$ and $\mathcal{N}$ the open subset of points $Q$ $\epsilon$ $M$,
such that there exists a neighborhood of every $Q$, where $\beta\neq0$.
Since, $\beta$ is a smooth function on $M$, then $\mathcal{W}\cup\mathcal{N}$ is
an open and dense subset of $M$. In $\mathcal{W}$ $\xi$ is principal. Furthermore, we consider $\mathcal{V}$, $\Omega$ open subsets of $\mathcal{N}$:
$$\mathcal{V}=\{Q\;\;\epsilon\;\;\mathcal{N}:\alpha=0\;\;in\;\;a\;\;neighborhood\;\;of\;\;Q\},$$
$$\Omega=\{Q\;\;\epsilon\;\;\mathcal{N}:\alpha\neq0\;\;in\;\;a\;\;neighborhood\;\;of\;\;Q\},$$
where $\mathcal{V}\cup\Omega$ is open and dense in the closure of $\mathcal{N}$.

\begin{proposition}
Let M be a real hypersurface in $M_{2}(c)$, equipped with
Lie $\mathbb{D}$-parallel structure Jacobi operator. Then,
$\mathcal{V}$ is empty.
\end{proposition}
\textbf{Proof:} Let $\{U,\varphi U,\xi\}$ be an orthonormal basis on $\mathcal{V}$. The following relations hold, because of Lemma 3.1
\begin{eqnarray}
AU=\gamma' U+\delta'\varphi U+\beta\xi,\hspace{20pt} A\varphi U=\delta'U+\mu'\varphi U,\hspace{20pt} A\xi=\beta U
\end{eqnarray}
\begin{eqnarray}
\nabla_{U}\xi=-\delta' U+\gamma'\varphi U,\hspace{20pt}
\nabla_{\varphi U}\xi=-\mu' U+\delta'\varphi U,\hspace{20pt}
\nabla_{\xi}\xi=\beta\varphi U,
\end{eqnarray}
\begin{eqnarray}
\nabla_{U}U=\kappa'_{1}\varphi U+\delta'\xi,\hspace{20pt}
\nabla_{\varphi U}U=\kappa'_{2}\varphi U+\mu'\xi,\hspace{20pt}
\nabla_{\xi}U=\kappa'_{3}\varphi U,
\end{eqnarray}
\begin{equation}
\nabla_{U}\varphi U=-\kappa'_{1}U-\gamma'\xi,\hspace{5pt}
\nabla_{\varphi U}\varphi U=-\kappa'_{2}U-\delta'\xi,\hspace{5pt}
\nabla_{\xi}\varphi U=-\kappa'_{3}U-\beta\xi,
\end{equation}
where $ \gamma',\delta',\mu',\kappa'_{1},\kappa'_{2},\kappa'_{3}$ are
smooth functions on $\mathcal{V}$.

From (2.6) for $X=U$ and $X=\varphi U$, taking into account (4.1), we obtain:
\begin{eqnarray}
l\varphi U=\frac{c}{4}\varphi U\hspace{20pt}lU=(\frac{c}{4}-\beta^{2})U.
\end{eqnarray}

Relation (1.1) , because of (4.2), (4.3) (4.4) and (4.5) implies:
\begin{eqnarray}
\delta'&=&0, \mbox{for\;\; $X=\varphi U$\;\;and\;\;$Y=\xi$}\\
(\mu'-\kappa'_{3})(\frac{c}{4}-\beta^{2})&=&0, \mbox{for\;\;$X=\varphi U$\;\;and\;\;$Y=\xi$}\\
\kappa'_{3}&=&\gamma', \mbox{for\;\;X=U\;\;and\;\;$Y=\xi$}.
\end{eqnarray}

On $\mathcal{V}$, relations (3.5)-(3.12), taking into account (4.6), become:
\begin{eqnarray}
\beta\kappa'_{1}+\mu'\kappa'_{3}+\frac{c}{4}&=&\gamma'\mu'+\gamma'\kappa'_{3}+\beta^{2}\\
\kappa'_{3}&=&3\mu'\\
(\varphi U)\beta&=&\beta\kappa'_{1}+\frac{c}{2}-2\gamma'\mu'\\
(\varphi U)\gamma'&=&\kappa'_{1}\gamma'+\beta\gamma'+2\beta\mu'-\mu'\kappa'_{1}.
\end{eqnarray}

Due to (4.7), we consider the open subsets $\mathcal{V}$:
$$\mathcal{V}_{1}=\{P\;\;\epsilon\;\;\mathcal{V}:\beta^{2}\neq\frac{c}{4}\;\;in\;\;a\;\;neighborhood\;\;of\;\;P\},$$
$$\mathcal{V}'_{1}=\{P\;\;\epsilon\;\;\mathcal{V}:\beta^{2}=\frac{c}{4}\;\;in\;\;a\;\;neighborhood\;\;of\;\;P\},$$
where $\mathcal{V}_{1}\cup\mathcal{V}'_{1}$ is open and dense in the closure of $\mathcal{V}$. So in $\mathcal{V}_{1}$ we obtain: $\mu'=\kappa'_{3}$.

On $\mathcal{V}_{1}$, because of relations (4.8) and (4.10), we obtain $\mu'=\kappa'_{3}=\gamma'=0$.

Relation (1.1), for $X=U$ and $Y=\varphi U$, due to (4.3), (4.4) and (4.5), implies: $\kappa'_{1}=0$. Substituting in (4.9) $\mu'=\kappa'_{3}=\gamma'=\kappa'_{1}=0$, leads to: $\beta^{2}=\frac{c}{4}$, which is impossible on $\mathcal{V}_{1}$. So $\mathcal{V}_{1}$ is empty and $\beta^{2}=\frac{c}{4}$ holds on $\mathcal{V}$.

On $\mathcal{V}$, because of (4.8) and (4.10), we have $\gamma'=\kappa'_{3}=3\mu'$. Substituting the last two relations in (4.9), implies: $\beta\kappa'_{1}=9\mu'^{2}$. Differentiation of $\beta^{2}=\frac{c}{4}$ with respect to $\varphi U$ and taking into account  (4.13), $\gamma'=3\mu'$ and $\beta\kappa'_{1}=9\mu'^{2}$, yields: $c=-6\mu'^{2}$, which is a contradiction because $\beta^{2}=\frac{c}{4}$. Hence, $\mathcal{V}=\emptyset$.
\qed\\

In what follows we work in $\Omega$.

By using (2.6), because of (3.1), we obtain:
\begin{eqnarray}
\hspace{10pt}
lU=(\frac{c}{4}+\alpha\gamma-\beta^{2})U+\alpha\delta\varphi
U\hspace{20pt}l\varphi U=\alpha\delta
U+(\alpha\mu+\frac{c}{4})\varphi U
\end{eqnarray}

Relation (1.1)  because of (3.2), (3.3) and (3.4) implies:
\begin{eqnarray}
\delta(\alpha\kappa_{3}+\frac{c}{4}-\beta^{2})&=&0 \;\;\mbox{for\;\;$X=U$\;\;and\;\;$Y=\xi$}\\
(\frac{c}{4}+\alpha\mu)(\kappa_{3}-\gamma)+\alpha\delta^{2}&=&0 \;\;\mbox{for\;\;$X=U$\;\;and\;\;$Y=\xi$}\\
(\frac{c}{4}+\alpha\gamma-\beta^{2})(\mu-\kappa_{3})-\alpha\delta^{2}&=&0\;\;\mbox{for\;\;$X=\varphi U$\;\;and\;\;$Y=\xi$}\\
\delta(\alpha\kappa_{3}+\frac{c}{4})&=&0\;\;\mbox{for\;\;$X=\varphi U$\;\;and\;\;$Y=\xi$}
\end{eqnarray}

Due to (4.17), we consider the open subsets $\Omega_{1}$ and $\Omega'_{1}$ of $\Omega$:
$$\Omega_{1}=\{Q\;\;\epsilon\;\;\Omega:\;\;\delta\neq0\;\;in\;\;a\;\;neighborhood\;\;of\;\;Q\},$$
$$\Omega'_{1}=\{Q\;\;\epsilon\;\;\Omega:\;\;\delta=0\;\;in\;\;a\;\;neighborhood\;\;of\;\;Q\},$$
where $\Omega_{1}\cup\Omega'_{1}$ is open and dense in the closure of $\Omega$.

In $\Omega_{1}$, from (4.14) and (4.17), we have: $\beta=0$, which is a contradiction, therefore $\Omega_{1}=\emptyset$. Thus we have: $\delta=0$ in $\Omega$ and relations from Lemma 3.1, (4.13), (4.15) and (4.16) become respectively:

\begin{eqnarray}
&&AU=\gamma U+\beta\xi,\hspace{20pt} A\varphi U=\mu\varphi U, \hspace{20pt} A\xi=\alpha\xi+\beta U\\
&&\nabla_{U}\xi=\gamma\varphi U,\hspace{20pt} \nabla_{\varphi U}\xi=-\mu U,\hspace{20pt} \nabla_{\xi}\xi=\beta\varphi U\\
&&\nabla_{U}U=\kappa_{1}\varphi U,\hspace{20pt}
\nabla_{\varphi U}U=\kappa_{2}\varphi U+\mu\xi,\hspace{20pt}
\nabla_{\xi}U=\kappa_{3}\varphi U,\\
&&\nabla_{U}\varphi U=-\kappa_{1}U-\gamma\xi,\hspace{5pt}
\nabla_{\varphi U}\varphi U=-\kappa_{2}U,\hspace{5pt}
\nabla_{\xi}\varphi U=-\kappa_{3}U-\beta\xi,\\
&&
lU=(\frac{c}{4}+\alpha\gamma-\beta^{2})U\hspace{20pt}l\varphi U=(\alpha\mu+\frac{c}{4})\varphi U\\
&&(\frac{c}{4}+\alpha\gamma-\beta^{2})(\mu-\kappa_{3})=0\\
&&(\frac{c}{4}+\alpha\mu)(\kappa_{3}-\gamma)=0.
\end{eqnarray}

Owing to (4.23), we consider the open subsets $\Omega_{2}$ and $\Omega'_{2}$ of $\Omega$:
$$\Omega_{2}=\{Q\;\;\epsilon\;\;\Omega:\;\;\mu\neq\kappa_{3}\;\;in\;\;a\;\;neighborhood\;\;of\;\;Q\},$$
$$\Omega'_{2}=\{Q\;\;\epsilon\;\;\Omega:\;\;\mu=\kappa_{3}\;\;in\;\;a\;\;neighborhood\;\;of\;\;Q\},$$
where $\Omega_{2}\cup\Omega'_{2}$ is open and dense in the closure of $\Omega$.

So in $\Omega_{2}$, we have:
\begin{eqnarray}
\gamma=\frac{\beta^{2}}{\alpha}-\frac{c}{4\alpha}\mbox{\;\;and\;\;(4.21)\;\;implies:\;\;$lU=0$}.
\end{eqnarray}

Due to (4.24) we consider $\Omega_{21}$ and $\Omega'_{21}$ the open subsets of $\Omega_{2}$:
$$\Omega_{21}=\{Q\;\;\epsilon\;\;\Omega_{2}:\;\;\mu=-\frac{c}{4\alpha}\;\;in\;\;a\;\;neighborhood\;\;of\;\;Q\},$$
$$\Omega'_{21}=\{Q\;\;\epsilon\;\;\Omega_{2}:\;\;\mu\neq-\frac{c}{4\alpha}\;\;in\;\;a\;\;neighborhood\;\;of\;\;Q\},$$
where $\Omega_{21}\cup\Omega'_{21}$ is open and dense in the closure of $\Omega_{2}$.
So in $\Omega_{21}$, (4.22) implies: $l\varphi U=0$ and because of Proposition 3.2 we obtain $\Omega_{21}=\emptyset$, thus in $\Omega_{2}$, we have $\mu\neq-\frac{c}{4\alpha}$ and as a result: $l\varphi U\neq0$. Furthermore, because of (4.24) $\kappa_{3}=\gamma$.

\begin{lemma}
Let M be a real hypersurface in $M_{2}(c)$, equipped with Lie $\mathbb{D}$-parallel structure Jacobi operator. Then $\Omega_{2}$ is empty.
\end{lemma}
\textbf{Proof:} 
On $\Omega_{2}$, relation (1.1) for $X=\varphi U$ and $Y=U$ owing to (4.20) and (4.21) implies: $\kappa_{2}l\varphi U=0$. So $\kappa_{2}=0$.
Due to the last, relations (3.8) and (3.12) imply:
\begin{eqnarray}
U\mu=\xi\mu=0.
\end{eqnarray}

Relation (1.1) for $X=U$ and $Y=\varphi U$ taking into account (4.20), (4.21), (4.22) and (4.25)  implies: $\mu=-\gamma$ and $\kappa_{1}=0$. Substitution in (3.6) of the relations which hold on $\Omega_{2}$ implies: $\gamma=0$ and this results in: $\beta^{2}=\frac{c}{4}$. Differentiation of the last with respect to $\varphi U$, because of (3.10) leads to $c=0$, which is a contradiction. This completes the proof of the present Lemma.
\qed\\

Summarizing, in $\Omega$ we have : $\delta=0$ and $\mu=\kappa_{3}$. Due to (4.24), we consider $\Omega_{3}$ and $\Omega'_{3}$ the open subsets of $\Omega$:
$$\Omega_{3}=\{Q\;\;\epsilon\;\;\Omega:\;\;\mu\neq-\frac{c}{4\alpha}\;\;in\;\;a\;\;neighborhood\;\;of\;\;Q\},$$
$$\Omega'_{3}=\{Q\;\;\epsilon\;\;\Omega:\;\;\mu=-\frac{c}{4\alpha}\;\;in\;\;a\;\;neighborhood\;\;of\;\;Q\},$$
where $\Omega_{3}\cup\Omega'_{3}$ is open and dense in the closure of $\Omega$. Since $\mu\neq-\frac{c}{4\alpha}$, due to (4.22) and (4.24) we obtain: $l\varphi U\neq0$ and $\gamma=\mu=\kappa_{3}$, in $\Omega_{3}$.
\begin{lemma}
Let M be a real hypersurface in $M_{2}(c)$, equipped with Lie $\mathbb{D}$-parallel structure Jacobi operator. Then $\Omega_{3}$ is empty.
\end{lemma}
\textbf{Proof:} 
 In $\Omega_{3}$ relations (3.6), (3.9), (3.10) and (3.11) become respectively:
\begin{eqnarray}
\alpha\gamma+\beta\kappa_{1}+\frac{c}{4}&=&\beta^{2}+\gamma^{2}\\
(\varphi U)\alpha&=&\alpha\beta-2\beta\gamma\\
(\varphi U)\beta&=&2\alpha\gamma+\beta\kappa_{1}+\frac{c}{2}-2\gamma^{2}\\
(\varphi U)\mu=(\varphi U)\gamma&=&3\beta\gamma.
\end{eqnarray}

Relation (1.1) for $X=Y=\varphi U$, because of (4.22), (4.28) and (4.30) implies: $\gamma(2\alpha-\gamma)=0$. Owing to the last relation, we consider  $\Omega_{31}$ and $\Omega'_{31}$ the open subsets of $\Omega_{3}$:
$$\Omega_{31}=\{Q\;\;\epsilon\;\;\Omega_{3}:\;\;\gamma\neq0\;\;in\;\;a\;\;neighborhood\;\;of\;\;Q\},$$
$$\Omega'_{31}=\{Q\;\;\epsilon\;\;\Omega_{3}:\;\;\gamma=0\;\;in\;\;a\;\;neighborhood\;\;of\;\;Q\},$$
where $\Omega_{31}\cup\Omega'_{31}$ is open and dense in the closure of $\Omega_{3}$.
In $\Omega_{31}$, $\gamma=2\alpha$. Differentiation of the latter with respect to $\varphi U$ and taking into  account (4.28) and (4.30) leads to: $\alpha\beta=0$, which is impossible. So $\Omega_{31}$ is empty.

Resuming on $\Omega_{3}$ we have: $\gamma=\mu=\kappa_{3}=0$ and relation (3.8) implies: $\kappa_{2}=0$. Relation (1.1) for $X=U$ and $Y=\varphi U$, because of (4.20), (4.21) and (4.22) yields: $\kappa_{1}=0$ and so relation (4.27) implies: $\beta^{2}=\frac{c}{4}$. Differentiation of the last along $\varphi U$ and because of (4.29) leads to $c=0$, which is a contradiction and this completes the proof of Lemma 4.3.
\qed\\

So on $\Omega$ the following relations hold:
$$\delta=0,\;\;\mu=\kappa_{3}=-\frac{c}{4\alpha},$$
and relation (4.22), because of the last one implies: $l\varphi U=0$.

Relation (1.1), for $X=U$ and $Y=\varphi U$, because of (4.20) and (4.21) implies: $\kappa_{1}lU=0$. So  $\kappa_{1}=0$, due to Proposition 3.2.

Due to the above relations, on $\Omega$ relations (3.9), (3.10) and (3.11) become respectively:
\begin{eqnarray}
(\varphi U)\alpha&=&\alpha\beta+\frac{c\beta}{2\alpha}\\
(\varphi U)\beta&=&\alpha\gamma+\frac{c\gamma}{2\alpha}+\frac{c}{4}\\
(\varphi U)\gamma&=&\beta\gamma-\frac{c\beta}{2\alpha}
\end{eqnarray}

Relation (1.1) for $X=\varphi U$ and $Y=U$ taking into account (4.20), (4.21) and (4.22) yields:
\begin{eqnarray}
(\varphi U)(\frac{c}{4}+\alpha\gamma-\beta^{2})&=&0\\
\kappa_{2}(\frac{c}{4}+\alpha\gamma-\beta^{2})&=&0\\
(\mu+\gamma)(\frac{c}{4}+\alpha\gamma-\beta^{2})&=&0.
\end{eqnarray}

Due to (4.35), we consider: $\Omega_{4}$ and $\Omega'_{4}$ the open subsets of $\Omega$:
$$\Omega_{4}=\{Q\;\;\epsilon\;\;\Omega:\;\;\frac{c}{4}+\alpha\gamma-\beta^{2}=0\;\;in\;\;a\;\;neighborhood\;\;of\;\;Q\},$$
$$\Omega'_{4}=\{Q\;\;\epsilon\;\;\Omega:\;\;\frac{c}{4}+\alpha\gamma-\beta^{2}\neq0\;\;in\;\;a\;\;neighborhood\;\;of\;\;Q\},$$
where $\Omega_{4}\cup\Omega'_{4}$ is open and dense in the closure of $\Omega$.
So in $\Omega_{4}$ we have: $\gamma=\frac{\beta^{2}}{\alpha}-\frac{c}{4\alpha}$ and because of (4.22) $lU=0$, which is impossible due to Proposition 3.2. Therefore, $\Omega_{4}=\emptyset$

So in $\Omega$ we have: $\kappa_{2}=0$ and relation (4.36) implies: $\mu=-\gamma$.

\begin{lemma}
Let M be a real hypersurface in $M_{2}(c)$, equipped with Lie $\mathbb{D}$-parallel structure Jacobi operator. Then $\Omega$ is empty.
\end{lemma}
\textbf{Proof:}  In $\Omega$ relations (3.8) and (3.12) yields: $U\alpha=\xi\alpha=0$.\\
Using the above relations, we obtain:
$$[U,\xi]\alpha=U\xi\alpha-\xi U\alpha=0,$$
$$[U,\xi]\alpha=(\nabla_{U}\xi-\nabla_{\xi}U)\alpha=\frac{c}{2\alpha}(\varphi U)\alpha.$$
Combining the last two relations and taking into account (4.31), we have: $c=-2\alpha^{2}$ and so $\mu=\frac{\alpha}{2}$ and $\gamma=-\frac{\alpha}{2}$. Relation (4.33), because of (4.31) and the last two relations imply: $\alpha=0$, which is impossible in $\Omega$. This completes the proof of Lemma 4.4.
\qed\\

We lead to the following due to Lemmas 4.1 and 4.4:
\begin{proposition}
Every real hypersurface in $M_{2}(c)$, equipped with Lie $\mathbb{D}$-parallel structure Jacobi operator is a Hopf hypersurface.
\end{proposition}
\section{Proof of Main Theorem}
Since M is a Hopf hypersurface, due to Theorem 2.1, (\cite{NR1}) , we have that $\alpha$ is a constant. We consider a unit vector field $Z$ $\epsilon$ $ker(\eta)$, such that $AZ=\lambda Z$, then $A\varphi Z=\nu\varphi Z$. Then $\{\xi, Z, \varphi Z\}$ is an orthonormal basis and the following relation holds on M, (Corollary 2.3, \cite{NR1}):
\begin{eqnarray}
\lambda\nu=\frac{\alpha}{2}(\lambda+\nu)+\frac{c}{4}
\end{eqnarray}
The relation (2.6) implies:
\begin{eqnarray}
lZ=(\frac{c}{4}+\alpha\lambda)Z\hspace{20pt}l\varphi Z=(\frac{c}{4}+\alpha\nu)\varphi Z
\end{eqnarray}
Relation (1.1) for $X=Z$ and $Y=\varphi Z$ and for $X=\varphi Z$ and $Y=Z$, because of (5.2) and taking the inner product of them with $\xi$ implies respectively:
\begin{eqnarray}
(\lambda+\nu)(\frac{c}{4}+\alpha\nu)&=&0\\
(\nu+\lambda)(\frac{c}{4}+\alpha\lambda)&=&0
\end{eqnarray}
I. Suppose that $\alpha\neq0$.

Due to  (5.3) and (5.4), we consider the open subset of $M$:
$$M_{1}=\{P\;\;\epsilon\;\;M:\lambda\neq-\nu\;\;in\;\;a\;\;neighborhood\;\;of\;\;P\}.$$
Because of (5.3) and (5.4) in $M_{1}$ we have $\lambda=\nu=-\frac{c}{4\alpha}$. In $M_{1}$ (5.1) yields $c=0$, which is a contradiction. Therefore, $M_{1}=\emptyset$.

Hence, in M we have: $\lambda=-\nu$. Substitution of the latter in (5.1) implies $c=-4\lambda^{2}$. From the last relation we conclude that: $c<0$ and $\lambda=constant$. The only hypersurface that we have in this case is of type B in $\mathbb{C}H^{2}$. Substituting the eigenvalues of this hypersurface in $\lambda=-\nu$ leads to a contradiction.

II. Suppose $\alpha=0$.

Relation (5.3) implies $\lambda=-\nu$ and so from relation (5.1) we obtain: $c=-4\lambda^{2}$. From the last two relations, we conclude that the only case which occurs is that of a real hypersurface in $\mathbb{C}H^{2}$ with three distinct constant eigenvalues. So it should be of type B in $\mathbb{C}H^{2}$, but for such hypersurface $\alpha$ can not vanish. So we lead to a contradiction and this completes the proof of our main theorem.

\end{document}